\documentstyle{article}
\newtheorem{definition}{Definition}[section]

\newtheorem{theorem}[definition]{Theorem}

\def \trait (#1) (#2) (#3){\vrule width #1pt height #2pt depth #3pt}
\def \qed{\hfill
        \trait (0.1) (6) (0)
        \trait (6) (0.1) (0)
        \kern-6pt
        \trait (6) (6) (-5.9)
        \trait (0.1) (6) (0)
\medskip}

\newcommand{\dt}{\,dt}

\newcommand{\e}{\varepsilon}

\newcommand{\FF}{{\cal F}}

\font\tenmsb=msbm10
\font\sevenmsb=msbm7
\font\fivemsb=msbm5
\newfam\msbfam
\textfont\msbfam=\tenmsb
\scriptfont\msbfam=\sevenmsb
\scriptscriptfont\msbfam=\fivemsb
\def\Bbb#1{{\fam\msbfam\relax#1}}

\def\rr{\Bbb R}

\newcommand{\NN}{{\Bbb N}}

\newcommand{\ZZ}{{\Bbb Z}}

\def\to{\rightarrow}

\def\H{{\cal H}^2}
\def\an{{\lambda_n}}
\def\xin{{x^i_n}}
\def\xip{{x^{i+1}_n}}
\def\xim{{x^{i-1}_n}}

\title{Discrete approximation of functionals with jumps and creases}

\author{ {\sc Andrea Braides}\\ SISSA Via Beirut 4 
  34014 Trieste, Italy}

\date{}

\begin{document}

\maketitle

\noindent Let $\H(0,L)$ denote the space of piecewise-$H^2$ functions on 
the interval $(0,L)$; 
if $u\in\H(0,L)$ then $u$, $u'$ and $u''$ are regarded 
as defined on the whole interval $(0,L)$, and $u$ and $u'$ are
piecewise-continuous functions. 
Let $S(u)$ be the set of the discontinuity points 
({\it jump points}) of the function 
$u\in\H(0,L)$ and, with a slight abuse of notation, denote 
by $S(u')$ the set of {\it crease points} of $u$ (i.e., those 
points where $u$ is continuous but $u'$ is discontinuous).
Let $\alpha,\beta>0$. The functional $\FF:\H(0,L)\to [0,+\infty)$
defined by
\begin{equation}\label{FF}
\FF(u)=\int_0^L |u''|^2\dt+\alpha\#(S(u'))+\beta\#(S(u))
\end{equation}
has been introduced by Blake and Zisserman 
(\cite{BZ}) to model some 
signal reconstruction problems, which can then be reduced  to
the solution of the minimum problems
\begin{equation}\label{pbbz}
m=\min\Bigl\{\FF(u)+\int_0^L|u-g|^2\dt:\ u\in\H(0,L)\Bigr\},
\end{equation}
where $g$ is the input (distorted) signal. The theoretical study 
of these problems has been 
performed by Coscia \cite{CO}, who interpreted the functional $\FF$
in the spirit of free-discontinuity problems as introduced
by De Giorgi and Ambrosio (see \cite{AFP}, \cite{BR}). 
The key point is to notice that if
$\alpha\le\beta\le 2\alpha$ then $\FF$ is lower semicontinuous 
with respect to the $L^1(0,L)$ convergence. At this point the
direct methods of the Calculus of Variations may be applied to 
obtain existence of the solutions to the abovementioned 
problems. A characterization of lower semicontinuity 
for general functionals on jumps and creases is given by Braides \cite{B}.

Here we show an approximation result by $\Gamma$-convergence of the 
functional above in the same spirit of that proved by Chambolle
\cite{Ch} for the Mumford-Shah functional (see \cite{MS}, \cite{MoS}, \cite{AFP})
$$
\int_0^L |u'|^2\dt + \gamma\#(S(u))
$$
defined on piecewise-$H^1$ functions. In that case, the approximating 
discrete functionals take the form
\begin{equation}\label{a1}
E_\e(u)=\sum_{x\in\e\ZZ\cap(0,L)}\e\Phi_\e\Bigl({u(x+\e)-u(x)\over\e}\Bigr),
\end{equation}
defined on discrete functions $u:\e\ZZ\cap(0,L)\to\rr$,
where $\Phi_\e$ are suitable functions, whose crucial property is to satisfy
$$
\lim_\e \Phi_\e(z)= z^2 \hbox{ on }\rr\qquad
\lim_\e \e\Phi_\e\Bigl({z\over\e}\Bigr)=\gamma
$$
($z\neq 0$). The model case is with $\Phi_\e(z)=\min\{z^2,\gamma/\e\}$.
The description of the behaviour of general
difference schemes of the form (\ref{a1}) with  minimal
hypotheses on $\Phi_\e$ can be found in \cite{BGe}.
In our case we can define suitable $\Psi_\e$ such that,
setting
$$
E_\e(u)=\sum_{x\in\e\ZZ\cap(0,L)}
\e\Psi_\e\Bigl({u(x+\e)+u(x-\e)-2u(x)\over\e^2}\Bigr)
$$
(i.e., using a discretization of $u''$ in place of the difference quotients)
these discrete energies converge to $\FF$ as $\e\to 0$
in the sense of $\Gamma$-convergence (see below). 
As a result, we immediately obtain approximate minimum problems
for (\ref{pbbz}) of the form
\begin{equation}\label{abz}
m_\e=\min\Bigl\{E_\e(u)+\sum_{x\in\e\ZZ\cap(0,L)}
\e|u(x)-g_\e(x)|^2:\ u:\e\ZZ\cap(0,L)\to \rr\Bigr\},
\end{equation}
where $g_\e$ are suitable discretizations of $g$.
The main 
requirement on $\Psi_\e$ is now that
$$
\lim_\e \Psi_\e(z)= z^2 \hbox{ on }\rr\qquad
\lim_\e \e\Psi_\e\Bigl({z\over\e}\Bigr)=\alpha,\qquad
\lim_\e 2\e\Psi_\e\Bigl({z\over\e^2}\Bigr)=\beta,\qquad
$$
thus highlighting an interesting double-scale effect.
Note in particular that the choice of 
$\Psi_\e(z)=\min\{z^2,\gamma/\e\}$
gives $\alpha=\gamma$ and $\beta=2\gamma$.

\bigskip
Before proceeding in the exact statement and proof of the result,
for notational convenience we replace the continuous 
small parameter $\e$
by a discrete parameter $\an$.
Let $L>0$ and $n\in\NN$. We set 
$$
\an={L\over n},\qquad 
\xin=i\an\quad(i=0,\ldots,n).
$$
Let $0<\alpha\le\beta\le 2\alpha$, $c_1, c_2>0$ and define
\begin{equation}\label{Psin}
\Psi_n(z)=\cases{\displaystyle z^2 & if $|z|\le {c_1/\sqrt{\an}}$
\cr\cr\displaystyle
{\alpha\over\an} & if ${c_1/\sqrt{\an}}<|z|\le
{c_2/\an\sqrt{\an}}$
\cr\cr\displaystyle
{\beta\over 2\an} & if ${c_2/\an\sqrt{\an}}<|z|$.}
\end{equation}
Note that the terms $c_1/\sqrt{\an}$ and $c_2/\an\sqrt{\an}$
may be replaced by any other pair of sequences $(T^1_n), 
(T^2_n)$ such that $1<\!<T^1_n<\!< 1/\an<\!< T^2_n<\!< 1/\an^2$
as $n\to+\infty$. Again, we have made this particular choice 
for notational convenience.

Let
$$
{\cal A}_n=\{u:\an\ZZ\cap[0,L]\to\rr\}
$$
and let
\begin{equation}\label{En}
E_n(u)=\sum_{i=1}^{n-1}\an\Psi_n
\Bigl({u(\xip)+u(\xim)-2u(\xin)\over\an^2}\Bigr)
\end{equation}
be defined for $u\in{\cal A}_n$.

We identify ${\cal A}_n$ with a subspace of $L^2(0,L)$
by regarding each function $u$ as defined by the value 
$u(\xin)$ on $[\xin,\xip)$. We then have the following
approximation result
by $\Gamma$-convergence. For a general introduction to $\Gamma$-convergence we refer
to \cite{GCB},
\cite{DM}, or \cite{BDF} Part II. For the application of $\Gamma$-convergence 
to the approximation of free-discontinuity problems we refer to \cite{BR}

\begin{theorem}
The functionals $E_n$ $\Gamma$-converge to $\FF$ with respect to the $L^1(0,L)$ 
convergence on bounded sets of $L^2(0,L)$. Namely, we have

{\rm(i)} if $(u_n)$ is bounded in $L^2(0,L)$ and $\sup_n E_n(u_n)<+\infty$ 
then $(u_n)$ is precompact in $L^1(0,L)$;

{\rm(ii)} if $(u_n)$ is bounded in $L^2(0,L)$,
 $u_n\to u$ in $L^1(0,L)$ and $\liminf_n E_n(u_n)<+\infty$ then 
$u\in\H(0,L)$ and $\FF(u)\le\liminf_n E_n(u_n)$;

{\rm(iii)} for all $u\in\H(0,L)$ there exist $u_n\in {\cal A}_n$
with $u_n\to u$ in $L^2(0,L)$ and $\FF(u)\ge\limsup_n E_n(u_n)$.

As a consequence of {\rm(i)--(iii)}, if $g\in L^2(0,L)$ then the minimum values
$$
m_n=\min\Bigl\{E_n(u)+\int_0^L |g-u|^2\dt: u\in {\cal A}_n\Bigr\}
$$
converge to the minimum value
$$
m=\min\Bigl\{\FF(u)+\int_0^L |g-u|^2\dt: u\in\H(0,L)\Bigr\}
$$
and from every sequence $(u_n)$ of solutions of the first
problems we can extract a converging subsequence to a solution
of the latter.
\end{theorem}

\begin{proof}
Statements (i) and (ii) will be proven by using the compactness
properties of $\FF$ and by comparing $E_\e$ with a suitable family
$(F_\e)$ $\Gamma$-converging to $\FF$. Let $u_n\in{\cal A}_n$ be
such that $\liminf_n E_n(u_n)<+\infty$. Upon extracting a subsequence,
we can suppose that this liminf is actually a limit,
so that in particular $\sup_n E_n(u_n)<+\infty$.
We will modify $u_n$ so as to obtain a comparison sequence $(v_n)$
in $\H(0,L)$. 

We first modify $u_n$ to obtain piecewise-affine functions $w_n$.
On $[\xim,\xin]$ the function $w_n$ is defined as
$$
w_n(x)= u_n(\xim)+(x-\xim)\Bigl({u_n(\xin)-u_n(\xim)\over\an}\Bigr).
$$
Let
$$
I_n=\{ i\in\{1,\ldots,n-1\}:\ u_n(\xip)+u_n(\xim)-2u_n(\xin)>c_1\an\sqrt{\an}\}.
$$
For all $i\in\{1,\ldots,n-1\}\setminus I_n$ we define $v_n$ on the
interval $(\xin-(\an/2),\xin+(\an/2))$ as the polynomial of degree $2$ satisfying
$$
v_n(\xin-{\an\over2})=w_n(\xin-{\an\over2})={u_n(\xin)+u_n(\xim)\over\an},
$$
$$
v_n(\xin+{\an\over2})=w_n(\xin+{\an\over2})={u_n(\xip)+
u_n(\xin)\over\an},
$$
$$
v'_n(\xin-{\an\over2})=w'_n(\xin-{\an\over2})={u_n(\xin)-u_n(\xim)\over\an},
$$
$$
v'_n(\xin+{\an\over2})=w'_n(\xin+{\an\over2})={u_n(\xip)-u_n(\xin)\over\an}.
$$
For such $v_n$ the constant second derivative
is given on $(\xin-(\an/2),\xin+(\an/2))$ by
$$
v''_n(y)= 
{1\over\an}(w'_n(\xin+)-w'_n(\xin-))
={u_n(\xip)+u_n(\xim)-2u_n(\xin)\over\an^2}.
$$
On the remaining part of $[0,l]$ we simply set
$$
v_n(x)=w_n(x).
$$
Note that $v_n\in\H(0,L)$ but $S(v_n)=\emptyset$. Moreover, 
$$
S(v'_n)=\{ \xin:\ i\in I_n\}.
$$
If $\xin\in S(v'_n)$ then 
\begin{equation}\label{aa1}
{u_n(\xip)+u_n(\xim)-2u_n(\xin)\over\an^2}
={1\over\an}(v'_n(\xin+)-v'_n(\xin-)).
\end{equation}
Finally, to complete the description of $v_n$, note that 
on the two intervals $(0,\an/2)$ and $(L-\an/2,L)$ 
$v_n$ is affine.

Let 
$$
F_n(v)=\int_0^L |v''|^2\dt
+\sum_{t\in S(v)}\varphi_n(v'(t+)-v'(t-))+\beta\#(S(v)).
$$
be defined on $\H(0,L)$, where 
$$
\varphi_n(z)=\cases{\alpha & if $|z|<c_2/\sqrt{\an}$
\cr\cr
\beta/2 & otherwise.}
$$
By construction, we immediately obtain that
$$ E_n(u_n)\ge F_n(v_n).$$
In fact, recalling (\ref{aa1}), we have
\begin{eqnarray*}
E_n(u_n)&=&\sum_{i\in\{1,\ldots,n-1\}\setminus I_n}
\an\Bigl|{u_n(\xip)+u_n(\xim)-2u_n(\xin)\over\an^2}\Bigr|^2
\\
&&+\alpha\,
\#\{i\in I_n:\ |u_n(\xip)+u_n(\xim)-2u_n(\xin)|<c_2\sqrt{\an}\}
\\
&&+{\beta\over 2}\,
\#\{i\in I_n:\ |u_n(\xip)+u_n(\xim)-2u_n(\xin)|\ge c_2\sqrt{\an}\}
\\
&=&\int_0^L|v_n''|^2\dt+\alpha\,
\#\{t\in S(v'_n):\ |v_n(t+)-v_n(t-)|<c_2/\sqrt{\an}\}
\\
&&+{\beta\over 2}\,
\#\{t\in S(v'_n):\ |v_n(t+)-v_n(t-)|\ge c_2/\sqrt{\an}\}
\\
&=&\int_0^L|v_n''|^2\dt+\sum_{t\in S(u)}
\varphi_n(v_n'(t+)-v_n'(t-))
\\
&=& F_n(v_n).
\end{eqnarray*}
It is easily checked that the $\Gamma$-limit of $F_n$ is $\FF$.

Note that $u_n-v_n$ tends to $0$ in $L^1(0,L)$. 
Since $F_n(v_n)\ge{\alpha\over\beta}\FF(v_n)$,
by the coerciveness properties of $\FF$ we obtain that $(v_n)$ is precompact
in $L^1(0,L)$
and each limit of a converging subsequence belongs
to $\H(0,L)$  (see \cite{CO}). 
Moreover, if $u_n\to u$ then also $v_n\to u$, and
by the $\Gamma$-convergence of $F_n$ to $\FF$ we obtain 
$$
\liminf_n E_n(u_n)\ge\liminf_n F_n(v_n)\ge \FF(u),
$$
and (ii) is proved.

To prove (iii), by a density argument it suffices to
show it for $u$ piecewise $C^2$ and with bounded second derivative.
Then, upon a slight piecewise-affine change of variable tending uniformly 
to the identity as $n\to+\infty$, we can reason as if
for all $n$ we have $S(u)\cup S(u')\subset \an\ZZ$.
We can then take $u_n(x)=u(x-)$ on $\an\ZZ$ 
(except at $0$, where we set $u_n(0)=u(0+)$);
i.e., we choose as $u_n$ in (iii) the piecewise-constant interpolation
of $u$. 
We have:

(a) if $\xin=t\in S(u)$ or $\xim=t\in S(u)$ then 
$$
{u_n(\xip)+u_n(\xim)-2u_n(\xin)\over\an^2}= 
{u(t+)-u(t-)+o(1)\over\an^2};
$$

(b) if $\xin\in S(u')$ then 
$$
{u(\xip)+u(\xim)-2u(\xin)\over\an^2}= 
{u'(\xin+)-u'(\xin-)+o(1)\over\an};
$$

(c) in all other cases, 
$$
{u(\xip)+u(\xim)-2u(\xin)\over\an^2}= 
{u''(\xin)+o(1)},
$$
with all the rests tending to $0$ uniformly as $n\to+\infty$.

Note that we have, in case (a),
$$
\Bigl|{u(\xip)+u(\xim)-2u(\xin)\over\an^2}|>\!> {c_2\over\an\sqrt{\an}},
$$
while in case (b) 
$$
{c_1\over\sqrt{\an}}<\!<
\Bigl|{u(\xip)+u(\xim)-2u(\xin)\over\an^2}|
<\!< {c_2\over\an\sqrt{\an}}
$$
as $n\to+\infty$.
Taking into account (a)--(c) and the remark above,
we immediately obtain
$$
\limsup_n E_n(u)\le \FF(u)
$$
so that (iii) is proved.

The final statement of the theorem follows easily from the 
well-known property of convergence of minima of $\Gamma$-limits.
In fact, from (ii) and the lower semicontinuity of
$u\mapsto\int_0^L|u-g|^2\dt$ with respect to the $L^1(0,L)$
convergence, we obtain $m\le\liminf_n m_n$. On the other hand, if $u$
is a solution of $m$, by (iii) we can find $(u_n)$ as in (iii)
so that 
$$
m=\lim_n\Bigl(E_n(u_n)+\int_0^L|u_n-g|^2\dt\Bigr)
\ge\limsup_n m_n.
$$
Finally, the pre-compactness property of minimizing sequences
follows from (i).
\end{proof}

\end{document}